\documentclass{svproc}
\usepackage{url}
\usepackage{amsfonts}   
\usepackage{amsmath}
\usepackage{bm}

\begin{document}
\mainmatter                        

\title{Spectral theory of matrix-sequences:
perspectives of the GLT analysis and beyond \thanks{Version \today}}

\titlerunning{GLT and beyond}

\author{Stefano Serra-Capizzano\inst{1,2}}

\authorrunning{S.~Serra-Capizzano}                   

\tocauthor{Stefano Serra-Capizzano}

\institute{Department of Science and High Technology,
Insubria University, Como, Italy\\
\email{s.serracapizzano@uninsubria.it}
\and
Department of Information Technology,
Uppsala University, Uppsala, Sweden\\
\email{serra@it.uu.se}}

\maketitle
\begin{abstract}
In recent years there has been a growing attention on distribution results in the sense of Weyl for the collective behavior of eigenvalues and singular values of matrix-sequences. Starting from the work of Szeg\"o regarding the case of Toeplitz matrix-sequences, there has been a wealth of associated results, which culminated in the works of Tilli and of Tyrtyshnikov, Zamarashkin, and of the author for preconditioned and non-preconditioned $r$-block $d$-level Toeplitz matrix-sequences with Lebesgue integrable generating functions. In the latter the use of matrix-valued linear positive operators and related Korovkin theories has been crucial. The subsequent steps induced by the analysis of preconditioning techniques and inspired by the rich world of the (pseudo) differential operators have been studies from the same perspective of matrix-sequences with hidden (asymptotic) structure, widely studied in the literature.
The widest generalization is represented by the notion of generalized locally Toeplitz matrix-sequences, which have inherent hidden structure and which include virtually any approximation via local numerical methods of (systems of) integral equations, partial and fractional differential equations, also with nonsmooth variable coefficients and irregular bounded/unbounded domains/manifolds.

In the current work, instead of focusing on a specific type of results, starting from the most recent advances on the topic, we describe shortly a series of open problems, challenges to be developed in the future.
\end{abstract}


\section{The GLT theory and around}\label{sec:intro}

	The idea of the generalized locally Toeplitz (GLT) matrix-sequences has been introduced in \cite{Serra03,Serra06}
as a generalization both of classical Toeplitz sequences and of variable coefficient (pseudo) differential operators \cite{hormander}. More specifically, for every sequence of the class, it has been demonstrated that it is possible to give a rigorous description
of the asymptotic spectrum in the Weyl sense \cite{BS,Ty96} in terms of a function (the symbol), that can be easily identified, via specific topological notions \cite{Tilli,Serra-acs} (see also \cite{equivGLT-topo,topo,glt-garoni,gacs} and references therein). The relevant notions, including that of approximating class of sequences \cite{Serra-acs} and of its recent generalizations \cite{gacs,blocking - irrational}, go back to the seminal paper by Tilli (and to his discussions with De Giorgi) \cite{Tilli}, where most of the important tools were also established, while related concepts were already present in the literature. See e.g. \cite{pre-LT}  for a classical result anticipating locally Toeplitz notions and \cite{ZamTy-ToeL1,Serra-taud1} for other related findings before the locally Toeplitz definitions.

The idea of the GLT symbol in fact is naturally connected with the notion of a symbol for (pseudo) differential operators - discrete and continuous - and, at the same time, it generalizes the notion of generating function for Toeplitz sequences, where for the latter the generating function is also the GLT symbol and it is identified through the Fourier coefficients, showing strong connections with the classical Fourier Analysis (see \cite{T-Chan SIREV}).

Below, for a given matrix-sequence of increasing sizes,  we report a definition in a general form of distributions of eigenvalues and singular values in the Weyl sense \cite{BS}.

\begin{definition}\label{def:sing and eig}
Let $\{A_n\}_n$ be a matrix sequence, where $A_n$ is of size $d_n$, $d_l <d_m$ when $l<m$, $l,m\in \mathbb{N}$, and let $\psi: D \subset \mathbb{R}^t \to \mathbb{C}^{r \times r}$ be a measurable function defined on a set $D$ with $0 < \mu_t(D) < \infty$, $\mu_t(\cdot)$ being the Lebesgue measure on 
$ \mathbb{R}^t$.
\begin{itemize}
    \item We say that $\{A_n\}_n$ has an \textbf{(asymptotic) singular value distribution} described by $\psi$, written as $\{A_n\}_n \sim_{\sigma} \psi$, if
    \begin{equation*}
        \lim_{n \to \infty} \frac{1}{d_n} \sum_{i=1}^{d_n} F(\sigma_i(A_n)) = \frac{1}{\mu_t(D)} \int_D \displaystyle\frac{\sum_{i=1}^{r} F(\sigma_i(\psi(\textbf{x})))}{r} d\textbf{x}, \quad \forall F \in C_c(\bbbr).
    \end{equation*}
    \item We say that $\{A_n\}_n$ has an \textbf{(asymptotic) spectral (eigenvalue) distribution} described by $\psi$, written as $\{A_n\}_n \sim_{\lambda} \psi$, if
    \begin{equation*}
        \lim_{n \to \infty} \frac{1}{d_n} \sum_{i=1}^{d_n} F(\lambda_i(A_n)) = \frac{1}{\mu_t(D)} \int_D \displaystyle\frac{\sum_{i=1}^{r} F(\lambda_i(\psi(\textbf{x})))}{r} d\textbf{x}, \quad \forall F \in C_c(\bbbc).
            \end{equation*}
            In this case, the function $\psi$ is referred to as the \textit{eigenvalue (or spectral) symbol} of $\{A_n\}_n$.
\end{itemize}
If $\psi$ describes both singular value and eigenvalue distribution of $\{A_n\}_n$, we write $\{A_n\}_n \sim_{\sigma, \lambda} \psi$.\\
 
\end{definition}

We now introduce the notion of the approximating class of sequences and a related key result \cite{Serra-acs}: both can be found in \cite{glt-garoni}, while the foundation of the notion can be traced back to the papers \cite{ZamTy-ToeL1,Tilli}. For a general topological study of this notion the reader is referred to \cite{topo,topo-G} and to the book \cite[Chapter 5]{glt-garoni}.

\begin{definition}\label{def:acs}
\textbf{(Approximating class of sequences)}
Let $\{A_n\}_n$ be a matrix-sequence and let $\{\{B_{n,m}\}_n\}_m$ be a class of matrix-sequences, with $A_n$ and $B_{n,m}$ of size $d_n$, $d_n$ strictly increasing sequence of positive integers.
We say that $\{\{B_{n,m}\}_n\}_m$ is an approximating class of sequences (a.c.s.) for $\{A_n\}_n$ if the following conditions are met:
for every $m$, there exist $n_m, c(m),\omega(m)$, such that, for every $n \geq n_m$,
\begin{equation*}
    A_n = B_{n,m} + R_{n,m} + N_{n,m},
\end{equation*}
with
\begin{equation*}
    {\rm rank } R_{n,m} \leq c(m) d_n, \quad \text{and} \quad \|N_{n,m}\| \leq \omega(m),
\end{equation*}
where $n_m$, $c(m)$, and $\omega(m)$ depend only on $m$ and
\begin{equation*}
    \lim_{m \to \infty} c(m) = \lim_{m \to \infty} \omega(m) = 0.
\end{equation*}
$
\{ \{ B_{n,m} \}_n \}_m \xrightarrow{\text{a.c.s. wrt } m} \{ A_n \}_n
$ denotes that $ \{ \{ B_{n,m} \}_n \}_m $ is an a.c.s. for $ \{ A_n \}_n $.
\end{definition}
The following theorem represents the expression of a related convergence theory and it is a powerful tool used, for example, in the construction of the various GLT $*$-algebras; see \cite{Serra03,Serra06,glt-garoni,glt-garoni-vol2,Barb rectangular,Barb 1D,Barb dD}.
\begin{theorem}\label{th:fundamental acs}
Let $\{A_n\}_n, \{B_{n,m}\}_n$, with $m,n \in \mathbb{N}$, be matrix-sequences and let $\psi, \psi_m : D \subset \mathbb{R}^d \to \mathbb{C}$ be measurable functions defined on a set $D$ with positive and finite Lebesgue measure. Suppose that
\begin{enumerate}
    \item $\{B_{n,m}\}_n \sim_{\sigma} \psi_m$ for every $m$;
    \item $\{\{B_{n,m}\}_n\}_m \xrightarrow{\text{a.c.s. wrt } m} \{A_n\}_n$;
    \item $\psi_m \to \psi$ in measure.
\end{enumerate}
Then
\begin{equation*}
    \{A_n\}_n \sim_{\sigma} \psi.
\end{equation*}
Moreover, if all the involved matrices are Hermitian, the first assumption is replaced by $ \{B_{n,m}\}_n \sim_{\lambda} \psi_m \; \text{for every } m,$ and the other two are left unchanged, then $\{A_n\}_n \sim_{\lambda} \psi.$
\end{theorem}

Regarding eigenvalue distribution and Theorem \ref{th:fundamental acs}, a useful generalization replaces the assumption on the Hermitian character of the matrices $A_n$ with proper Schatten $p$ norm conditions on $A_n-A_n^*$ for $n$ large enough, with $X^*$ indicating the transpose conjugate of a matrix $X$. The reader is referred to \cite{GoSe} for this type of results in the context of eigenvalues of compact perturbations of Jacobi matrix-sequences and to \cite{Barb NonH} for refined results and applications to approximated PDEs; see also Section \ref{subsec:glt_algebra}, Axiom \textbf{GLT 5.}, and \cite{Barb dD}, for a connection with the GLT theory.

There exist several matrix-sequences with explicit or hidden (asymptotic) structure widely studied in the literature as reported in \cite{hidden}. Among them we have the \( d \)-level \( r \)-block Toeplitz matrix-sequences, which have an explicit structure and which enjoy Szeg\"o like relations as in Definition \ref{def:sing and eig}.
They have been studied deeply in the last century both for their mathematical beauty \cite{Tilli-ToeL1} and for their pervasive use in applications (see e.g. \cite{Toe-use1,Toe-use2} and references therein). A wide generalization is represented by the notion of generalized locally Toeplitz (GLT) matrix-sequences \cite{Tilli,Serra03,Serra06}, which have inherent hidden structure and whose construction is heavily based on Definition \ref{def:sing and eig}, Definition \ref{def:acs}, and Theorem \ref{th:fundamental acs}. 

Indeed, for every $r,d\ge 1$ the $r$-block $d$-level GLT class has nice $*$-algebra features and indeed it has been proven
that it is stable under linear combinations, products, and inversion when the sequence which is inverted
shows a sparsely vanishing symbol (sparsely vanishing symbol $=$ a symbol whose minimal singular value vanishes at most in a set
of zero Lebesgue measure). Furthermore, the GLT $*$-algebras virtually include any approximation of partial differential equations (PDEs), fractional differential equations (FDEs), integro-differential equations (IDEs) by local methods  (Finite Difference, Finite Element, Isogeometric Analysis  etc). The use of the GLT theory has been proven effective in studying the zeros of orthogonal polynomials \cite{GoSe,Kuij}, in spectral operator theory \cite{bianchi,bianchi-S,holm-serra,GoUMI,koro-infinity-1,koro-infinity-2}, in fast eigenvalue solvers for structured matrices \cite{sven-erik-asymp-exp-prec,BGS,BGS-2,eig expansion Laplacian,sven-erik-asymp-exp1}, in a PDE/FDE/IDE numerical setting for the efficient iterative solution of large linear systems \cite{BeSe,Bene Rognes,BogoFDE1,BogoFDE2,DGMSS,DGMSScol,our-MC,garoni-IgA-multigrid,fract-derivatives,fract-derivatives2,DoroNS,Ivo,eig expansion Laplacian,sven-erik-asymp-exp1,curl-curl IgA GLT,NS Schur GLT,fract-derivatives3,immersed NLAA-2}, including (systems of) PDEs/FDEs also with nonsmooth variable coefficients and irregular bounded domains/manifolds (see \cite{glt-garoni,glt-garoni-vol2,Barb 1D,Barb dD} and references therein).

The seminal papers on the considered spectral theory are \cite{T-Chan SIREV,Ty96,Tilli,Serra03,Serra06}.

The GLT theory is organized in books and revue papers \cite{glt-garoni,glt-garoni-vol2,glt-garoni-CIME,Barb rectangular,Barb 1D,Barb dD}. We mention also other references contain applications and theoretical results, where either the GLT machinery has been used, or new related techniques are introduced for the asymptotic eigenvalue analysis in a non-Hermitian/non-normal setting \cite{IE,Barb NonH,BGS,Ivo,flipped1,flipped PDE,flipped wave control,Ngondiep,Alec Beta,Alec cluster}. The new branch on momentary symbols \cite{mom-Toe,mom-GLT,dyad-new2} and on fast eigensolvers can be found in \cite{BogoFDE1,BGS,BGS-2,BogoFDE2}. The considered fast eigensolvers called `matrix-less' represent a fruitful connections between the asymptotic expansions studies in operator theory by Bogoya, B\"ottcher, Grudsky and the algorithmic proposals in numerical analysis by Ekstr\"om, Furci, Garoni, and the author.

\section*{Challenges}

There are also various possible directions in which the GLT theory could be expanded or used and which could form a program for future investigations. As possible items, we could  include symbol based preconditioning, multigrid, linear positive operators (LPOs) and Korovkin theory, spectral detection of branches, fast matrix-less computation of eigenvalues, stability issues, and new challenges such as the GLT use in tensors, stochastic, machine learning algorithms.

In the remainder of the work we stress also the impact and the further potential of the theory with special attention to new tools and to new directions as those based on symmetrization tricks \cite{flipped2,flipped1,flipped-BCs,flipping GLT1,flipping GLT2,flipped PDE,flipped wave control}, on the extra-dimensional approach \cite{draft-T,curl-curl IgA GLT,gacs,blocking - pre-prequel}, and on blocking structures/operations \cite{blocking - prequel,blocking - irrational}.

\section{Spectral tools}\label{sec:spectral}

In this section, we introduce the essential tools for the spectral analysis of the matrices under consideration, using the block multi-level GLT matrix sequence framework. The case of scalar values, corresponding to the non-block setting, has been extensively detailed in \cite{glt-garoni,glt-garoni-vol2}, while the block setting, associated with matrix-valued symbols, is thoroughly examined in \cite{Barb 1D,Barb dD}; see also \cite{Barb rectangular} for the rectangular setting.

\subsection{Notation and terminology}
\textbf{Matrices and matrix-sequences}. Given a square matrix $A \in \mathbb{C}^{m \times m}$, we denote by $A^*$ its conjugate transpose and by $A^\dagger$ the Moore–Penrose pseudoinverse of $A$. Recall that $A^\dagger = A^{-1}$ whenever $A$ is invertible. The singular values and eigenvalues of $A$ are denoted respectively by $ \sigma_1(A), \cdots, \sigma_m(A) $ and $\lambda_1(A), \cdots, \lambda_m(A) $.\\
Regarding matrix norms, $\|\cdot\|$ refers to the spectral norm, and for $1 \leq p \leq \infty$, the notation $\|\cdot\|_p$ stands for the Schatten $p$-norm defined as the $p$-norm of the vector of singular values. Note that the Schatten $\infty$-norm, which is equal to the largest singular value, coincides with the spectral norm $\|\cdot\|$. The Schatten $1$-norm is often referred to as the trace-norm, being the sum of the singular values, while the Schatten $2$-norm coincides with the Frobenius norm. Schatten $p$-norms, as important special cases of unitarily invariant norms, are treated in detail in a wonderful book by Bhatia \cite{Bhatia1997}.\\
Finally, the expression \textbf{matrix-sequence} refers to any sequence of the form $\{A_n\}_n$, where $A_n$ is a square matrix of size $d_n$ with $d_n$ strictly increasing so that $d_n \to \infty$ as $n \to \infty$. A \textbf{r-block matrix-sequence}, or simply a matrix-sequence if $r$ can be deduced from context is a special $\{A_n\}_n$ in which the size of $A_n$ is $d_n = r\varphi_n$, with $r \geq 1 \in \mathbb{N}$ fixed and $\varphi_n \in \mathbb{N}$ strictly increasing.
\subsection{Multi-index notation}
In order to deal effectively with multilevel structures, it is very convenient to use multi-indices, which are vectors of the form $\textbf{i} = (i_1, \cdots, i_d) \in \mathbb{Z}^d$. The related notation is listed below and it is taken from \cite{glt-garoni-vol2}, which borrowed it from the seminal paper by Tyrtyshnikov on spectral distribution via matrix-theoretic tools \cite{Ty96}.
\begin{itemize}
    \item $\bm{0}, \bm{1}, \bm{2}, \dots$ are vectors of all zeroes, ones, twos, etc.
    \item $\textbf{h} \leq \textbf{k}$ means that $h_r \leq k_r$ for all $r = 1, \cdots, d$. In general, relations between multi-indices are evaluated componentwise.
    \item Operations between multi-indices, such as addition, subtraction, multiplication, and division, are also performed componentwise.
    \item The multi-index interval $[\textbf{h}, \textbf{k}]$ is the set $\{\textbf{j} \in \mathbb{Z}^d : \textbf{h} \leq \textbf{j} \leq \textbf{k}\}$. We always assume that the elements in an interval $[\textbf{h}, \textbf{k}]$ are ordered in the standard lexicographic manner
  \[
\Biggr[ \cdots
\biggr[\Bigr[\bigl(j_1, \cdots, j_d\bigl)\Bigr]
_{j_d = h_d, \cdots, k_d}\biggr]
_{j_{d-1} = h_{d-1}, \cdots, k_{d-1}}
\cdots \Biggr]_{j_1 = h_1, \cdots, k_1}.
\]
    \item $\textbf{j} = \textbf{h}, \cdots, \textbf{k}$ means that $\textbf{j}$ varies from $\textbf{h}$ to $\textbf{k}$, always following the lexicographic ordering.
    \item $\textbf{m} \to \infty$ means that $\min(\textbf{m}) = \min_{j=1, \cdots, d} m_j \to \infty$.
    \item The product of all the components of $\textbf{m}$ is denoted as $\nu(\textbf{m}) := \prod_{j=1}^{d} m_j$.
\end{itemize}
A \textbf{multilevel matrix-sequence} is a matrix-sequence $\{A_{\textbf{n}} \}_n$ such that $n$ varies in some infinite subset of $\mathbb{N}$, $\textbf{n} = \textbf{n}(n)$ is a multi-index in $\mathbb{N}^d$ depending on $n$, and $\textbf{n} \to \infty$ when $n \to \infty$. This is typical of many approximations of differential operators in $d$ dimensions.\\
\textbf{Measurability.} All the terminology from measure theory, such as “measurable set”, “measurable function”, “a.e.”, etc., refers to the Lebesgue measure in $ \mathbb{R}^t $, denoted with $\mu_t $. A matrix-valued function $f : D \subseteq \mathbb{R}^t \to \mathbb{C}^{r \times r} $ is said to be measurable (resp., continuous,
Riemann-integrable, in $ L^p(D) $, etc.) if all its components $ f_{\alpha\beta} : D \to \mathbb{C} $,  $ \alpha, \beta = 1, \dots, r $, are measurable (resp., continuous, Riemann-integrable, in $ L^p(D) $, etc.). If $ f_m, f : D \subseteq \mathbb{R}^t \to \mathbb{C}^{r \times r} $ are measurable, we say that $ f_m $ converges to $ f $ in measure (resp., a.e., in $ L^p(D) $, etc.) if $ (f_m)_{\alpha\beta} $ converges to
$ f_{\alpha\beta} $ in measure (resp., a.e., in $ L^p(D) $, etc.) for all $ \alpha, \beta = 1, \cdots, r $. If $A \in \mathbb{C}^{m \times m}$, the singular values and eigenvalues of $A$ are denoted by $\sigma_1(A), \dots, \sigma_m(A)$ and $\lambda_1(A), \dots, \lambda_m(A)$, respectively. The set of eigenvalues (i.e., the spectrum) of $A$ is denoted by $\Lambda(A)$.

\subsection{ Singular Value and Eigenvalue Distributions of a Matrix-Sequence}

With reference to Definition \ref{def:sing and eig}, it is useful to recall that the informal meaning behind the spectral distribution definition is the following: if $\psi$  is continuous, then a suitable ordering of the eigenvalues $\{\lambda_j(A_n)\}_{j=1,\cdots,d_n}$,
 assigned in correspondence with an equispaced grid on $D$, reconstructs approximately the $r$ surfaces $\textbf{x} \to \lambda_i(\psi(\textbf{x})),  i =1,...,r.$ For instance, in the simplest case where $t =1$ and $D =[a,b]$,
$d_n = nr$, the eigenvalues of $A_n$ are approximately equal - up to few potential outliers to  $\lambda_i (\psi(x_j)),$ where
 \begin{equation*}
  x_j = a + j \frac{(b-a)}{n}, \quad j =1,\cdots,n, \quad i =1,\cdots,r.
 \end{equation*}
If $t =2$ and $D = [a_1,b_1] \times [a_2,b_2]$,
$d_n = n^2r$, the eigenvalues of $A_n$ are approximately equal, again up to a few potential outliers, to $\lambda_i (\psi(x_{j_1},y_{j_2})), i =1,\cdots,r,$ where
 \begin{align*}
 x_{j_1} &= a_1 + j_1 \frac{(b_1 - a_1)}{n}, \qquad j_1 =1,\cdots,n,\\
 y_{j_2} &= a_2 + j_2 \frac{(b_2 - a_2)}{n},\qquad j_2 =1,\cdots,n,
 \end{align*}
If the considered structure is two-level, then the subscript is $\bm{n}=(n_1,n_2)$ and $d_n = n_1 n_2 r$. Furthermore, for $t \geq 3$, a similar reasoning applies.\\
Finally, we report an observation that is often useful in theoretical derivations.
\begin{remark}\label{rem: range}
The relation $\{A_n\}_n \sim_\lambda f$ and $\Lambda(A_n) \subseteq S$ for all $n$ imply that the range of $f$ is a subset of the closure $\bar{S}$ of
$S$. In particular, $\{A_n\}_n \sim_\lambda f$ and $A_n$ positive definite for all $n$ imply that $f$ is non-negative definite almost everywhere, simply nonnegative almost everywhere if $r=1$; see e.g. \cite{Barb NonH} for its use in the spectral analysis of approximated PDEs and \cite{GM-GLT} for the spectral study of matrix-sequences stemming from geometric means. The same applies when a multilevel matrix sequence $\{A_{\bm{n}}\}_{\bm{n}}$ is considered and similar statements hold when singular values are taken into account.
\end{remark}

With regard to the notion  of approximating class of sequences, we notice that Definition \ref{def:acs} and Theorem \ref{th:fundamental acs} both hold with obvious changes, when the involved matrix-sequences show a multilevel structure. In that case $n$ is replaced by $\bm{n}$ uniformly in $A_{n}, B_{n,m}, d_{n}$.

\subsection{Matrix-Sequences with Explicit or Hidden (Asymptotic) Structure} \label{subsec:matrix_structures}
In this subsection, we introduce three types of (asymptotic) matrix structures that serve as the fundamental building blocks of the GLT \(*\)-algebras and in fact they are the algebra generators. Specifically, for any positive integers \( d \) and \( r \), we consider the set of \( d \)-level \( r \)-block GLT matrix-sequences. This set forms a \(*\)-algebra of matrix-sequences, which is both maximal and isometrically equivalent to the maximal \(*\)-algebra of \( 2d \)-variate \( r \times r \) matrix-valued measurable functions (with respect to the Lebesgue measure) that are naturally defined over $ [0,1]^d \times [-\pi,\pi]^d; $ see \cite{glt-garoni,glt-garoni-vol2,Barb 1D,Barb dD} and references therein.\\
The reduced version of the GLT asymptotic structure \cite{Serra06,reduced} plays a crucial role in approximating integro-differential operators, including their fractional versions, particularly when defined over general (non-Cartesian) domains. This concept was initially introduced in \cite{Serra03,Serra06} and later extensively developed in \cite{reduced}, where GLT symbols are again defined as measurable functions over $ \Omega \times [-\pi,\pi]^d,$ with \( \Omega \) being Peano-Jordan measurable and contained within \( [0,1]^d \). Additionally, the reduced versions also form maximal \(*\)-algebras that are isometrically equivalent to their corresponding maximal \(*\)-algebras of measurable functions.\\
These GLT \(*\)-algebras provide a rich framework of hidden (asymptotic) structures, built upon two fundamental classes of explicit algebraic structures: \( d \)-level \( r \)-block Toeplitz matrix-sequences and sampling diagonal matrix-sequences (discussed in Sections \ref{sec:multi} and \ref{sec:Block}), along with asymptotic structures described by zero-distributed matrix-sequences (see Section \ref{sec:zero}). Notably, the latter class serves an analogous role to compact operators in relation to bounded linear operators, forming a two-sided ideal of matrix-sequences within any of the GLT \(*\)-algebras.

\subsection{Zero-distributed sequences}\label{sec:zero}
 Zero-distributed sequences are defined as matrix sequences $\{A_n\}_n$ such that $\{A_n\}_n \sim_\sigma 0$. Note that, for any $r \geq 1$, $\{A_n\}_n \sim_\sigma 0$ is equivalent to
$\{A_n\}_n \sim_\sigma O_r$, where $O_r$ is the $r \times r$ zero matrix. The following theorem, taken from \cite{Serra-taud2,glt-garoni}, provides a useful characterization for detecting this type of sequence.
\begin{theorem}
Let $\{A_n\}_n$ be a matrix-sequence, with $A_n$ of size $d_n$ and let $p\in [1,\infty]$, with $\|X\|_p$ being the Schatten $p$-norm of $X$, that is the $l^p$ norm of the vector its singular values. Let $\|\cdot\|=\|\cdot\|_\infty$ be the spectral norm. With
the natural convention $1/\infty = 0$, we have
\begin{itemize}
    \item $\{A_n\}_n \sim_\sigma 0$ if and only if $A_n = R_n + N_n$ with ${\operatorname{rank}(R_n)}/{d_n} \to 0$ and $ ||N_n|| \to 0$ as $n \to \infty$;
    \item $\{A_n\}_n \sim_\sigma 0$ if there exists $p \in [1,\infty]$ such that
    \begin{equation*}
        \frac{\|A_n\|_p}{(d_n)^{1/p}} \to 0 \quad \text{as} \quad n \to \infty.
    \end{equation*}
\end{itemize}
\end{theorem}
Notice that the same definition can be given and corresponding result (with obvious changes) holds, when the involved matrix-sequences show a multilevel structure.
In that case $n$ is replaced by $\bm{n}$ uniformly in $A_{n}, N_{n}, R_{n}, d_{n}$.

Further results connecting Schatten $p$-norms, unitarily invariant norms (see \cite{Bhatia1997}), variational characterization, and Toeplitz/Hankel structures can be found in \cite{Serra-Tilli2}.

Going back the definition of zero-distributed matrix-sequences, taking into consideration the equivalence $\{A_n\}_n \sim_\sigma 0$ if and only if $A_n = R_n + N_n$ with ${\operatorname{rank}(R_n)}/{d_n} \to 0$ and $ ||N_n|| \to 0$ as $n \to \infty$, we observe that the zero-distributed matrix-sequences coincide with the important class of asymptotical low-rank matrix-sequences. Furthermore, they include the numerical approximation of integral operators, those modeling blurring operators in signal processing and imaging, dense sampling matrix-sequences, Hankel matrix-sequences generated by a function in their maximal generality, structured matrix-sequences stemming from the extra-dimensional approach and from the blocking techniques; see e.g. \cite{IE,FT,draft-T,blocking - pre-prequel,blocking - irrational,blocking - num,blocking - prequel} and references therein and \cite{SSS} for a specific application in finance. 

\subsection{ Multilevel block Toeplitz matrices}\label{sec:multi}

Fix $r,d\ge 1$ positive integers. Given \(\textbf{n} \in \mathbb{N}^d \), a matrix of the form
\begin{equation*}
[A_{\textbf{i}-\textbf{j}}]_{\textbf{i,j=1}}^{\textbf{n}} \in \mathbb{C}^{\nu(\textbf{n})r \times \nu(\textbf{n})r},
\end{equation*}
with blocks \( A_\textbf{k} \in \mathbb{C}^{r \times r} \), \( \textbf{k} \in \{-(\textbf{n}-1), \dots, \textbf{n}-1\} \), is called a multilevel block Toeplitz matrix, or, more precisely, a \( d \)-level \( r \)-block Toeplitz matrix.\\
Given a matrix-valued function \( f : [-\pi,\pi]^d \to \mathbb{C}^{r \times r} \) belonging to \( L^1([- \pi, \pi]^d) \), the \( \mathbf{n} \)-th Toeplitz matrix associated with \( f \) is defined as
\begin{equation*}
T_\mathbf{n}(f):=[\hat{f}_{\mathbf{i-j}}]_{\mathbf{i,j=1}}^{\mathbf{n}} \in \mathbb{C}^{\nu(\mathbf{n})r \times \nu(\mathbf{n})r},
\end{equation*}
where
\begin{equation*}
\hat{f}_\mathbf{k} = \frac{1}{(2\pi)^d} \int_{[-\pi,\pi]^d} f(\bm{\theta}) e^{-\hat{\iota} (\mathbf{k}, \bm{\theta)}} d\bm{\theta} \in \mathbb{C}^{r \times r}, \quad \mathbf{k} \in \mathbb{Z}^d,
\end{equation*}
are the Fourier coefficients of \( f \), in which \( \hat{\iota} \) denotes the imaginary unit, the integrals are computed componentwise, and \( (\mathbf{k}, \bm{\theta}) = k_1\theta_1 + \cdots + k_d\theta_d \). Equivalently, \( T_\mathbf{n}(f) \) can be expressed as
\begin{equation*}
T_\mathbf{n}(f) =
\sum_{|j_1|<n_1} \cdots \sum_{|j_d|<n_d} [J_{n_1}^{(j_1)}
\otimes \cdots \otimes J^{(j_d)}_{n_d}] \otimes \hat{f}(j_1, \cdots, j_d),
\end{equation*}
where \( \otimes \) denotes the Kronecker tensor product between matrices, and \( J^{(l)}_{m} \) is the matrix of order \( m \) whose \( (i,j) \) entry equals 1 if \( i - j = l \) and zero otherwise.\\
The family \( \{T_{\mathbf{n}}(f)\}_{\mathbf{n}\in\mathbb{N}^d} \) is the family of (multilevel block) Toeplitz matrices associated with \( f \), which is called the generating function.

\subsection{Multilevel block diagonal sampling matrices}\label{sec:Block}
 Given $r,d \geq 1$ positive integers, given $\mathbf{n} \in \mathbb{N}^d$ and a function $a : [0,1]^d \to \mathbb{C}^{r \times r}$, we define the multilevel block diagonal sampling matrix $D_{\mathbf{n}}(a)$ as the block diagonal matrix
\begin{equation*}
    D_{\mathbf{n}}(a) = \underset{\mathbf{i=1,\dots,n}}{\operatorname{diag}} a \left( \mathbf{\frac{i}{n}} \right) \in \mathbb{C}^{\nu(\mathbf{n})r \times \nu(\mathbf{n})r}.
\end{equation*}

\subsection{The $*$-slgebra of $d$-Level $r$-block GLT Matrix-Sequences}\label{subsec:glt_algebra}

Let $r,d \geq 1$ be fixed positive integers. A $d$-level $r$-block GLT sequence, or simply a GLT sequence if we do not need to specify either $r$ or $d$, is a special $d$-level $r$-block matrix-sequence  equipped with a measurable function $ \kappa :[0,1]^d \times[-\pi,\pi]^d \to \mathbb{C}^{r\times r}, \quad d \geq 1,$
called the symbol. The symbol is essentially unique, in the sense that if $\kappa, \xi$ are two symbols of the same
GLT sequence, then $\kappa = \xi$ a.e. We write \( \{A_n\}_n \sim_{\mathrm{GLT}} \kappa \) to denote that $\{A_n\}_n$ is a GLT
sequence with symbol $\kappa$.\\
It can be proven that the set of multilevel block GLT sequences is the $*$-algebra
generated by the three classes of sequences defined in Section \ref{sec:zero}, Section \ref{sec:multi}, Section \ref{sec:Block}: zero-distributed,
multilevel block Toeplitz, and block diagonal sampling matrix sequences. The GLT $*$-algebras satisfy several algebraic and topological properties that are treated in detail in \cite{glt-garoni,glt-garoni-vol2,Barb 1D,Barb dD}.

As mentioned at the beginning of the present work, there are several related results and applications to the approximation via local numerical methods of (systems of) PDEs/FDEs also with nonsmooth variable coefficients and irregular bounded domains/manifolds; for the spectral analysis in the case of $d$-level $p$-block GLT asymptotic structures see \cite{Barb 1D,Barb dD,DNS,DoroNS,NS DG Dumb}, for GLT based fast numerical solvers also on non-Cartesian domains and variable coefficients see
\cite{DGMSS,DGMSScol,our-MC,fract-derivatives2,curl-div IgA GLT,curl-curl IgA GLT}, while in \cite{immersed NLAA-1,immersed NLAA-2} general domains and trimmed geometries are considered, the review paper \cite{hughes} contains an engineering perspective, and the work \cite{GLT-PDE-manifold} includes the GLT analysis in the case of PDEs on manifolds.

Here, we focus on the main operative properties listed below that represent a complete characterization of GLT sequences, equivalent to the full constructive definition.

\subsection*{GLT Axioms}

\begin{itemize}
    \item \textbf{GLT 1.} If $\{A_{\bm{n}}\}_{\bm{n}} \sim_{\mathrm{GLT}} \kappa$ then $\{A_{\bm{n}}\}_{\bm{n}} \sim_\sigma \kappa$ in the sense of Definition \ref{def:sing and eig}, with $t = 2d$ and $D = [0, 1]^d \times [-\pi, \pi]^d$. Moreover, if each $A_{\bm{n}}$ is Hermitian, then $\{A_{\bm{n}}\}_{\bm{n}} \sim_\lambda \kappa$, again in the sense of Definition \ref{def:sing and eig} with $t = 2d$.
    \item \textbf{GLT 2.} We have
    \begin{itemize}
        \item $\{T_{\bm{n}}(f)\}_{\bm{n}} \sim_{\mathrm{GLT}} \kappa(\bm{x}, \bm{\theta}) = f(\bm{\theta)}$ if $f : [-\pi, \pi]^d \to \mathbb{C}^{r \times r}$ is in $L^1([- \pi, \pi]^d)$;
        \item $\{D_{\bm{n}}(a)\}_{\bm{n}} \sim_{\mathrm{GLT}} \kappa(\bm{x}, \bm{\theta}) = a(\bm{x})$ if $a : [0, 1]^d \to \mathbb{C}^{r \times r}$ is Riemann-integrable;
        \item $\{Z_{\bm{n}}\}_{\bm{n}} \sim_{\mathrm{GLT}} \kappa(\bm{x}, \bm{\theta}) = O_r$ if and only if $\{Z_{\bm{n}}\}_{\bm{n}} \sim_\sigma 0$.
    \end{itemize}

    \item \textbf{GLT 3.} If $\{A_{\bm{n}}\}_{\bm{n}} \sim_{\mathrm{GLT}} \kappa$ and $\{B_{\bm{n}}\}_{\bm{n}} \sim_{\mathrm{GLT}} \xi$, then:
    \begin{itemize}
        \item $\{A_{\bm{n}}^*\}_{\bm{n}} \sim_{\mathrm{GLT}} \kappa^*$;
        \item $\{\alpha A_{\bm{n}} + \beta B_{\bm{n}}\}_{\bm{n}} \sim_{\mathrm{GLT}} \alpha \kappa + \beta \xi$ for all $\alpha, \beta \in \mathbb{C}$;
        \item $\{A_{\bm{n}} B_{\bm{n}}\}_{\bm{n}} \sim_{\mathrm{GLT}} \kappa \xi$;
        \item $\{A_{\bm{n}}^\dagger\}_{\bm{n}} \sim_{\mathrm{GLT}} \kappa^{-1}$, provided that $\kappa$ is invertible almost everywhere.
    \end{itemize}
    \item \textbf{GLT 4.} \(\{A_{\bm{n}}\}_{\bm{n}} \sim_{\mathrm{GLT}} \kappa\) if and only if there exist \(\{B_{{\bm{n}},j}\}_{\bm{n}} \sim_{\text{GLT}} \kappa_j\) such that \(\{\{B_{{\bm{n}},j}\}_{\bm{n}}\}_j \xrightarrow{\text{a.c.s.wrt }j} \{A_{\bm{n}}\}_{\bm{n}}\) and \(\kappa_j \rightarrow \kappa\) in measure.
    \item \textbf{GLT 5.} If $\{A_{\bm{n}}\}_n \sim_{\mathrm{GLT}} \kappa$ and $A_{\bm{n}} = X_{\bm{n}} + Y_{\bm{n}}$, where
\begin{itemize}
    \item every $X_{\bm{n}}$ is Hermitian,
    \item $ ||X_{\bm{n}}||, \, ||Y_{\bm{{n}}}|| \leq C$ for some constant $C$ independent of $\bm{{n}}$,
    \item ${\nu}(\bm{n})^{-1} \|Y_{\bm{n}}\|_1 \to 0$,
\end{itemize}
then $\{A_{\bm{n}}\}_n \sim_\lambda \kappa$.
   \item \textbf{GLT 6.}  If $\{A_{\bm{n}}\}_n \sim_{\mathrm{GLT}} \kappa$ and each $A_{\bm{n}}$ is Hermitian, then $\{f(A_{\bm{n}})\}_n \sim_{\mathrm{GLT}} f(\kappa)$ for every continuous function $f : \mathbb{C} \to \mathbb{C}$.

\end{itemize}

Note that, by \textbf{GLT 1}, it is always possible to obtain the singular value distribution from the GLT symbol, while the eigenvalue distribution can only be deduced either if the involved matrices are Hermitian or the related matrix-sequence is quasi-Hermitian in the sense of Axiom \textbf{GLT 5}; see also \cite{Barb NonH,GoSe}.

Furthermore, Axiom \textbf{GLT 2}, part 3, is quite important because it means that all asymptotical low-rank matrix-sequences are GLT with zero symbol
\cite{IE,FT,draft-T,blocking - pre-prequel,blocking - irrational,blocking - num,blocking - prequel,SSS}.

\section{Challenges, conjectures, new directions}

The challenges indicated at the end of Section \ref{sec:intro} represent open directions for future researches. We describe them more explicitly in terms of ideas in the subsequent list below. They represent a vast research program, which is subdivided into 8 research projects: [RP 1.), [RP 2.), [RP 3.), [RP 4.), [RP 5.), [RP 6.), [RP 7.), [RP 8.).

\begin{description}
\item[RP 1. Stochastic GLT theory]
	The main target is to revise the whole GLT theory in the stochastic setting, so allowing to treat crucial applications modeled by stochastic PDEs (SPDEs). Despite the fact that the international community working in the area of SPDEs is relatively young, in the last years there has been an outstanding development of this field in the direction of both mathematical theory and applications. This has been also confirmed by important international recognitions starting from the Fields Medal awarded to Martin Hairer.
	
The first step is the extension of the definition of the approximating class of sequences in a stochastic version. This is already highly nontrivial given the presence of asymptotical ‘low rank’ terms.

A possibility is as follows.

\begin{definition}\label{def: S-acs}
\textbf{(Stochastic approximating class of sequences)}
Let $\{A_n\}_n$ be a matrix-sequence and let $\{\{B_{n,m}\}_n\}_m$ be a class of matrix-sequences, with $A_n$ and $B_{n,m}$ of size $d_n$, $d_n$ strictly increasing sequence of positive integers.
We say that $\{\{B_{n,m}\}_n\}_m$ is a stochastic approximating class of sequences (S-a.c.s.) for $\{A_n\}_n$ if the following conditions are met:
for every $m$, there exists $n_m,c(m),s(m),\omega(m)\ge 0$ such that, for every $n \geq n_m$,
\begin{equation*}
    A_n = B_{n,m} + S_{n,m} + R_{n,m} + N_{n,m},
\end{equation*}
with
\begin{eqnarray*}
{\rm prob}\left\{{\rm rank }(R_{n,m})/d_n\leq c(m)\right\} & > &1-1/m, \\
{\rm prob}\left\{ \|N_{n,m}\|\leq \omega(m)\right\} & > &1-1/m, \\
{\rm prob}\left\{S_{n,m}\neq O_{M_{d_n}}\right\} & \le & s(m),
\end{eqnarray*}
 where $n_m$, $c(m)$, $s(m)$, and $\omega(m)$ depend only on $m$ and
\begin{equation*}
    \lim_{m \to \infty} c(m) =  \lim_{m \to \infty} s(m) = \lim_{m \to \infty} \omega(m) = 0.
\end{equation*}
In short we write $\{ \{ B_{n,m} \}_n \}_m \xrightarrow{\text{S-a.c.s. wrt } m} \{ A_n \}_n$, for denoting that the class $ \{ \{ B_{n,m} \}_n \}_m $ is a S-a.c.s. for $ \{ A_n \}_n $.
\end{definition}

The second step is to adapt \cite[Theorem 5.4, Theorem 5.5, Theorem 5.6, Corollary 5.1, Corollary 5.2]{glt-garoni} to the stochastic setting where the a.c.s approximation is replaced by its stochastic counterpart, both in the eigenvalue and singular value sense. In this direction, the definitions of distribution have to be given in stochastic terms, i.e. $\sim_{S,\sigma},\sim_{S,\lambda}$, but, for that matter, the literature already furnishes few tools from the world of random matrices.

The following could be a formulation of the stochastic counterpart of Theorem \ref{th:fundamental acs} and could be employed in the construction of the corresponding stochastic $r$-block $d$-level GLT $*$-algebras, $d,r\ge 1$, mimicking the steps performed in \cite{Serra03,Serra06,glt-garoni,glt-garoni-vol2,Barb 1D,Barb dD} for their non-stochastic versions.

\begin{theorem}\label{th:fundamental S-acs}
Let $\{A_n\}_n, \{B_{n,m}\}_n$, with $m,n \in \mathbb{N}$, be matrix-sequences and let $\psi, \psi_m : D \subset \mathbb{R}^d \to \mathbb{C}$ be measurable functions defined on a set $D$ with positive and finite Lebesgue measure. Suppose that
\begin{enumerate}
    \item $\{B_{n,m}\}_n \sim_{S,\sigma} \psi_m$ for every $m$;
    \item $\{\{B_{n,m}\}_n\}_m \xrightarrow{\text{S-a.c.s. wrt } m} \{A_n\}_n$;
    \item $\psi_m \to \psi$ in measure.
\end{enumerate}
Then
\begin{equation*}
    \{A_n\}_n \sim_{S,\sigma} \psi.
\end{equation*}
Moreover, if all the involved matrices are Hermitian, the first assumption is replaced by $ \{B_{n,m}\}_n \sim_{S,\lambda} \psi_m \; \text{for every } m,$ and the other two are left unchanged, then $\{A_n\}_n \sim_{S,\lambda} \psi$.
\end{theorem}

Once step 1 and step 2 are performed, there is a whole freeway of potential results.
\begin{itemize}
\item The whole GLT theory in a stochastic setting, with all the variations of $d$-level $r$-block cases, with $d,r\ge 1$;
\item	 Then a really remarkable set of concrete applications can be treated, taking inspiration from the world of modeling stemming from SPDEs.
\end{itemize}
\item[RP 2. Automatic symbol computation] 	There are hundreds of papers showing that virtually any discretization of PDEs - on graded meshes by using any meaningful approximation method (Finite Differences, Finite Volumes, Finite Elements, Isogeometric Analysis, Discontinuous Galerkin and many more) - leads to $d$-level $r$-block GLT sequences, potentially of reduced type (see \cite[pp. 398-399]{Serra03}, \cite[Section 3.1.4]{Serra06}, \cite{reduced}), with appropriate $r=r({\cal F},d)$ and $d$, with $d$ dimensionality of the domain where the continuous operator is defined, and $r$  proper block size depending on $d$ and on the specific approximation formula $\cal F$.
Ratnani had the idea of a symbolic automatic computation of the GLT symbol (see \cite[Chapter 11, item 3]{glt-garoni} and \cite[Chapter 8, item 2]{glt-garoni-vol2}) where the initial proposal is mentioned). Then a still basic development is proposed in \cite{Sarathkumar}. However, we are far from a working general procedure: this task is challenging and it should be tried on a specific set of approximation procedures combined with specific applications modeled by PDEs, FDEs, IEs.

Connections between banded GLT structures, generalized Toeplitz graph sequences, zero distribution of variable coefficient orthogonal polynomials, and approximated PDEs by local methods can be found in \cite{BeSe,Kuij,graph-GLT-1,graph-GLT-2}: these special types of structures are easier to treat from the point of view of the identification of the GLT symbol, thanks to the finite Fourier expansion in the Fourier variables.

\item[RP 3. Machine Learning in the GLT theory] 	If a fixed GLT matrix-sequence is given with invertible symbol a.e. and a proposed preconditiong matrix-sequence is constructed having the same GLT symbol, then the difference matrix-sequence is still GLT (by the $*$-algebra structure of the GLT matrix-sequences) and its symbol is zero. Under very mild assumptions, the latter means that the preconditioned matrix-sequence has symbol $1$, which implies that the spectral approximation holds in a distributional sense.
However, the presence of $o(d_n)$ outliers is possible, which means that the preconditioned matrix-sequence can be written as $I_{rd_n}$ plus a zero-distributed matrix-sequence. A challenge, which would open a whole research direction, is to combine the GLT framework with machine learning techniques that learn the low rank structure of the zero-distributed component. While the distributional convergence remains unchanged, identifying and exploiting the rank structure via data‐driven methods could accelerate the convergence of such “data-informed” or “data-augmented” preconditioners.

\item[RP 4. Toward tensor-valued generating functions and GLT symbols] 	In a recent arXiv contribution \cite{toe-quaternion}, the authors study several spectral and numerical features of standard Toeplitz matrices and matrix-sequences, where the generating function takes values into the algebraic structure of quaternions. Beside the applicative side, this research line is of interest because it opens the door to several potential results. We divide them into 4 subitems.
\begin{itemize}
\item The work \cite{toe-quaternion} is not complete since the considered generating functions are essentially bounded: we could investigate the complete picture with Lebesgue integrable generating functions, with also the singular value distribution as in the seminal works by Tilli \cite{Tilli-ToeL1} and Tyrtyshnikov, Zamarashkin \cite{ZamTy-ToeL1} in the complex-valued case (even rectangular matrix-valued \cite{Tilli-ToeL1}).
\item The first subitem would complete the unilevel quaternion setting. However, the $d$-level $r\times s$ block quaternion Toeplitz setting is the next generalization, which corresponds to the case of a $d$-variate $r\times s$ matrix-valued generating function; see \cite{Barb rectangular} and references therein.
\item Both previous items could be complemented by adding variations in the physical domains. Again applications of interest have to be investigated, together with a complete $d$-level, $r\times s$ block quaternion GLT theory.
\item The three subitems described before all stem from the quaternion variation in the values of the generating functions. This invites for more variations, e.g. when the generating functions tale values in complex $m$-tensors or even in quaternion complex $m$-tensors, in one or several variables. The whole GLT theory in a tensor setting could be developed and this is valuable because, when considering the solution of stochastic SPDE, a way for obtaining fast algorithm is to compute it and to store it in a $3$-tensor. Its fast computation depends on the fact that the associated matrix/tensor-sequence is zero-distributed: the latter represents already a link between stochastic and tensor GLT theories.
\end{itemize}

\item[RP 5. Connections with the infinite dimensional setting]
Here we report few open questions in \cite[Section 3.11]{Serra06}, which are still open after 20 years.

As clear from the literature \cite{glt-garoni,glt-garoni-vol2,Barb 1D,Barb dD}, the structure of the definition of GLT sequences needs the notion of Toeplitz generated by a symbol as building block and we have observed that there is no difference depending on the nature of the symbol in the deduced results: the symbol can be in one variable or in many variables, it can be scalar-valued or (rectangular) matrix-valued. Therefore one may ask wether it is of interest building a GLT theory where the symbols of the Toeplitz structures are e.g. operator-valued.

A second question concerns the link between GLT theory and H\"ormander (pseudo) differential operators and between the general spectral results for approximations of PDEs and the spectrum of the associated continuous operators. It is clear that we have some technical difficulty in generalizing a distribution formula  in the case of the continuous operator, for instance in the case of an infinite countable unbounded spectrum. Under the latter hypothesis, to give a meaning to left hand-side of a distribution formula we should have something for replacing the matrix-size $d_n$; moreover, since $F$ has bounded support and in the important case of continuous operators having a discrete spectrum with $\infty$ as unique accumulation point (see e.g. the classical differential operators considered in {\rm \cite{question:kac}}), the summation will contain only a finite number of terms. A possibility can be as follows: fix $t$ positive value and consider
\begin{equation*}\label{distribution:continuous}
{1\over \#\{\lambda_n :\ |\lambda_n|<t\}}\sum_{|\lambda_n|<t}F(\lambda_n/t)
\end{equation*}
or, for second order operators,
\begin{equation*}\label{distribution:continuous-2}
{1\over t^{d/2}}\sum_{\lambda_n<t}F(\lambda_n/t)
\end{equation*}
since  $\#\{\lambda_n :\ |\lambda_n|<t\}= \#\{\lambda_n :\ \lambda_n<t\}\sim t^{d/2}$, where $F$ is a usual test function and the $\lambda_n$s are the eigenvalues of the considered differential operator with the right boundary conditions. The latter are both potential analogs of the standard spectral distribution in the Weyl sense,  and then the next step is to understand how the above quantity behaves as $t$ tends to infinity: in this direction, see the generalizations of Weyl's formula in {\rm \cite{Fle1,Fle2,Gil,Gru,Ivrii,See}}, the beautiful question by Kac {\rm \cite{question:kac}} (where, at page 4 setting $\lambda=t/v$, $v$ constant in $(0,1]$, (\ref{distribution:continuous}) occurs with $F$ being $vF_0(v\cdot)$ with $F_0(\cdot)$ being the characteristic function of $[0,1]$), the beautiful answer 25 years later in {\rm \cite{answer:kac}}, and the discussion in {\rm \cite{holm-serra}} and in \cite[Section 3.1]{Serra06}.

The connection emphasized above between discrete and continuous has many applications, including very popular ones in control theory; see for example the paper \cite{zuazua}. Recent works on this matter are \cite{bianchi-S,eig expansion Laplacian,bianchi,stile-bianchi}, but a general procedure for deducing the eigenvalue distribution of a continuous problem from (a sequence of) its numerical approximations (of increasing precision) is definitely open.

\item[RP 6. LPOs, Korovkin theory, statistical convergence, and GLT theory]\ \\
The Toeplitz and GLT matrix-sequences with quaternion-valued symbols are interesting also in connection with the theory of LPOs, because we could also consider the extension of the Korovkin theory in that setting, as done in the standard case in \cite{Serra-koro1,Serra-koro2,Serra-koro-lama,koro-infinity-1,koro-infinity-2}. The whole machinery would be based on a proper partial ordering and to related notion of (matrix-valued) LPO, which is essential in the Korovkin theory. We remind that the notion of matrix-valued LPO has an importance that goes beyond the Korovkin theory.
Indeed the distributional results for matrix-sequences can be studied, not only from the viewpoint of the GLT theory, but also from the viewpoint of matrix-valued LPOs and in that case we have also very useful pointwise localization results on the spectra, which can be exploited in the context of structure preserving optimal preconditioning strategies.
As further items in this direction we can mention the following two items.
\begin{itemize}
\item the extension of the outlier free approximations to variable coefficient elliptic operators, in the spirit of \cite{outl-free},  where the constant coefficient case is treated in detail. The main tools could be the LPO notion and the dyadic decomposition techniques already used in the relevant literature (see e.g. \cite{dyad-1,dyad-2,dyad-3,dyad-new1,dyad-new2}).
\item The challenge of using the clustering theory obtained using the matrix Korovkin results in the context of the statistical convergence \cite{stat-conv1,stat-conv2,stat-conv3,stat-conv4,stat-conv5}, which represents a whole sub-branch of the analytic side of the approximation theory.
\end{itemize}

\item[RP 7. Maximality GLT Conjecture]
This item can be metaphorically described as follows: are the GLT sisters the only other maximal $*$-algebras in the space of matrix-sequences of size 
$r\nu(\textbf{n})$?

As proven by Barbarino \cite{equivGLT-topo}, the $r$-block $d$-level GLT class is a maximal $*$-algebra in the space of sequences with sizes $\nu(\textbf{n})r\times \nu(\textbf{n})r$, being isometrically equivalent to $r\times r$ matrix-valued Lebesgue measurable functions in $[0,1]^d \times [-\pi,\pi]^d$ in Lebegue measure. The basic question is: with $r$ and $d$ fixed, is it true that, up to extra-dimensional compressions \cite{blocking - pre-prequel}, the only other maximal algebras are sisters of the corresponding $r$-block $d$-level GLT $*$-algebra? Here for GLT sister we mean the following. Given a fixed sequence  $\{U_\textbf{n}\}$ of unitary matrices and given $r,d\ge 1$, the new sister is obtained as the set of all matrix-sequences formed as $\{U_\textbf{n} X_\textbf{n} U_\textbf{n}^*\}$, where there exists $\textbf{n}_0$ and $U_\textbf{n}^* U_\textbf{n}=I$, $I$ identity of size $r\nu(\textbf{n}) =r\cdot n_1 \cdot n_2\cdots n_d$ for $\textbf{n}\ge \textbf{n}_0 $ and $\{X_\textbf{n}\}$ being any $r$-block $d$-level GLT matrix-sequence.
We notice that the question holds in the case of reduced GLT $*$-algebras as well.

We observe that it would be sufficient to start with the case $d=r=1$. The rest is technical.

\item[RP 8. Surfing on permuted GLT $*$-algebras using the blocking idea]\ \\
This item can be metaphorically described as surfing among GLT $*$-algebras and still maintaining a limit distribution, so defining an envelope containing all the transformed and reduced/compressed GLT $*$-algebras. For dealing with such a type of problems it is useful to remark that the extra-dimensional approach is an important tool. The initial proposal is given by Tyrtyshnikov in \cite{draft-T}, while few theoretical techniques are developed in \cite{curl-curl IgA GLT,blocking - pre-prequel,blocking - prequel,blocking - num,gacs}: in particular in \cite{gacs} a generalization of the a.c.s. notion is given, called generalized approximating class of sequences (g.a.c.s.).

We observe that Theorem 3.2 in \cite{blocking - irrational}  is important also for its interpretation in terms of $r$-block $d$-level GLT $*$-algebras (with $d=1$). In fact, by \cite{gacs}, we know that the g.a.c.s. limit as $m$ tends to infinity of  $\left\{\{B_{n,m}\}_{n}\right\}_m$ is $\{A_n\}_n $, the g.a.c.s. topology that can be described via the standard a.c.s. topology, after making all the matrices of the same size. Furthermore, by the main result in \cite{blocking - prequel}, $B_{n,m}$ can be chosen similar by permutations to a Toeplitz matrix $T_{n,m}$ generated by matrix-valued symbol, whose matrix-size $r_m$ depends on $m$. Hence, by \cite{Barb 1D}, $\{T_{n,m}\}_{n}$ is a $r_m$-block GLT matrix-sequence so that $\{A_n\}_n $ is obtained as the limit over $m$ of matrix-sequences which belong to permuted $r_m$-block GLT $*$-algebras. However, for any $r$, in general the matrix-sequence $\{A_n\}_n$ does not belong to the $r$-block GLT $*$-algebra.

Finally, the surfing problem can be considered also in a $d$-level setting for $d>1$. However this extension is technical: the true relevant question is represented by the case of $d=1$.

\end{description}

We end the section by observing that there exist further challenges that we omit to discuss in detail, like for instance a richer theoretical analysis concerning Toeplitz and GLT momentary symbols, especially when applied to the design of fast matrix-less eigensolvers \cite{mom-Toe,mom-GLT,dyad-new2,BGS-2}.

\section{Conclusions}

We have given a general essential overview of the GLT theory and of its main building blocks, like to the topological notion of a.c.s. convergence. This has been employed for describing several challenges and open research directions, which could be considered as inter-connected research projects.
We remind that in 2006 paper \cite{Serra06} had the same flavour, containing several potential new directions, including that of the reduced GLT matrix-sequences which has been completed by Barbarino \cite{reduced} in 2022. However, some of the questions are still open, such as that in item [RP 5.), and the present contribution can be regarded as an update and an enrichment of the program in \cite{Serra06}.

Regarding the various chapters, the research project in [RP 1.) is very appealing, but it is challenging because of technical difficulties. Regarding the foundation of the a.c.s. convergence, the tools that could be used can be found in the paper by Tao and Vu \cite{TaoVu-RandMatr}. For the eigenvalue distribution in a stochastic sense, we refer to \cite[Section 2.5]{book-random-matrices} is of interest for the the symmetric setting and in general \cite[Chapter 2]{book-random-matrices}, plus \cite[Chapter 4]{book-random-matrices} contains important notions and tools in a more general context, e.g. for the singular value distribution.

Research projects [RP 2.) and [RP 3.) are also challenging, but there is so much to do that the two projects can be modulated and adapted at the level of a young researcher.

The research project in [RP 4.) is the most adequate for an independent study even by a PhD student. It is very modular with several steps of increasing difficulty.

Research projects [RP 5.), [RP 6.), [RP 7.), [RP 8.)  have a more analytical flavour, but also the practical impact could be quite strong, especially for items [RP 5.) and [RP 6.).

Overall the given research projects represent a program of several years and must include the work of several researchers.

\hrule

\medskip
\begin{thebibliography}{99}
	%

	\def\BIT{{\it BIT\ }}
\def\BITNM{{\it BIT Numer.\ Math.\ }}
\def\CMA{{\it Comput.\ Math.\ Appl.\ }}
\def\cal{{\it Calcolo\ }}
\def\JCP{{\it J.\ Comput.\ Phys.\ }}
\def\CMAME{{\it Comput.\ Methods Appl.\ Mech.\ Engrg.\ }}
\def\JAT{{\it J.\ Approx.\ Theory\ }}
\def\LAA{{\it Linear Algebra Appl.\ }}
\def\LMA{{\it Linear\ Multilin.\ Algebra\ }}
\def\MC{{\it Math. Comp.\ }}
\def\MMMAS{{\it Math.\ Models Methods Appl.\ Sci.\ }}
\def\NLAA{{\it Numer.\ Linear Algebra Appl.\ }}
\def\NM{{\it Numer.\ Math.\ }}
\def\OM{{\it Oper.\ Matrices\ }}
\def\SIAMJMAA{{\it SIAM J.\ Matrix\ Anal.\ Appl.\ }}
\def\SIAMJNA{{\it SIAM J.\ Numer.\ Anal.\ }}
\def\SIAMJSC{{\it SIAM J.\ Sci.\ Comput.\ }}



\bibitem{graph-GLT-1} A. Adriani, D. Bianchi, S. Serra-Capizzano. Asymptotic spectra of large (grid) graphs with a uniform local structure (Part I): Theory. {\it Milan J. Math.} \textbf{88} (2020), no. 2, 409–454.

\bibitem{graph-GLT-2} A. Adriani, D. Bianchi, P. Ferrari, S. Serra-Capizzano. Asymptotic spectra of large (grid) graphs with a uniform local structure, Part II: Numerical applications. {\it J. Comput. Appl. Math.} \textbf{437} (2024), Paper No. 115461.


\bibitem{blocking - irrational}
A. Adriani, A.J.A. Schiavoni-Piazza, S. Serra-Capizzano.
Blocking structures, g.a.c.s. approximation, and distributions.
{\it Bol. Soc. Mat. Mex. (3)}  \textbf{31-2} (2025), Paper No. 41.

\bibitem{gacs}
A. Adriani, A.J.A. Schiavoni-Piazza, S. Serra-Capizzano, C. Tablino-Possio.
Revisiting the notion of approximating class of sequences for handling approximated PDEs on moving or unbounded domains. {\it Electron. Trans. Numer. Anal.} (2025), in press.

\bibitem{GLT-PDE-manifold}
A. Adriani, M. Semplice, S. Serra-Capizzano. {Generalized Locally Toeplitz matrix-sequences and approximated PDEs on submanifolds: the flat case.} {\it  Linear Multilin. Algebra} \textbf{73-9} (2025), 1902–1924.


\bibitem{sven-erik-asymp-exp-prec}
 F. Ahmad, E.S. Al-Aidarous, D. Abdullah Alrehaili, S.-E. Ekstr\"om, I. Furci, S. Serra-Capizzano.
 Are the eigenvalues of preconditioned banded symmetric Toeplitz matrices  known in almost closed form?
{\it  Numer. Alg.} \textbf{78-3} (2018), 867-893.

\bibitem{IE} A.S. Al-Fhaid, S. Serra-Capizzano, D. Sesana, M.Z. Ullah. Singular-value (and eigenvalue) distribution and Krylov preconditioning of sequences of sampling matrices approximating integral operators.
{\it  Numer. Linear Algebra Appl.} \textbf{21-6} (2014), 722-743.


\bibitem{book-random-matrices}
G.W. Anderson, A. Guionnet, O. Zeitouni. \emph{An introduction to random matrices}.
Cambridge Studies in Advanced Mathematics, 118. Cambridge University Press, Cambridge, 2010.


\bibitem{Arico-Donatelli-Serra}
{A. Aric\`o, M. Donatelli, S. Serra-Capizzano}.
{V-cycle optimal convergence for certain (multilevel) structured linear systems}.
\SIAMJMAA \textbf{26} (2004), 186-214.

\bibitem{stat-conv2} J. Banaś, M. Mursaleen. {\em Sequence spaces and measures of noncompactness with applications to differential and integral equations}. Springer, New Delhi, 2014.

\bibitem{blocking - num}
N. Barakitis, M. Donatelli, S. Ferri, V. Loi, S. Serra-Capizzano, R.L. Sormani.
Block structures, approximation, and preconditioning. {\it Numer. Alg.} (2025),  online 8-7-25.

\bibitem{blocking - pre-prequel}
N. Barakitis, P. Ferrari, I. Furci, S. Serra-Capizzano.
An extradimensional approach for distributional results: the case of $2\times 2$ block Toeplitz structures.
{\it Springer Proc. Math. Stat. on Mathematical Modeling with Modern Applications} \textbf{497} (2025), 61-78;  online 9-7-25. 

\bibitem{dyad-new2}
N. Barakitis, V. Loi, S. Serra-Capizzano. A note on eigenvalues and singular values of variable Toeplitz matrices and matrix-sequences with application to variable two-step BDF approximations to parabolic equations. {\it Proc. International Conference on Spectral and Approximation Theory}, (Kochin - Kerala - (India)  November 27th--30th,  2023), {\em Springer book series ``Trends in Mathematics"} (2025), 37-67; online 27-7-25-


\bibitem{equivGLT-topo}
G. Barbarino. Equivalence between GLT sequences and measurable functions. {\it Linear Algebra Appl.} \textbf{529}  (2017), 397–412.

\bibitem{reduced}
G. Barbarino. A systematic approach to reduced GLT. {\it  BIT} \textbf{62-3} (2022), 681-743.

\bibitem{topo} G. Barbarino, C. Garoni. From convergence in measure to convergence of matrix-sequences through concave functions and singular values. {\it Electron. J. Linear Algebra} \textbf{32} (2017), 500–513.

\bibitem{Barb rectangular} G. Barbarino, C. Garoni, M. Mazza, S. Serra-Capizzano. Rectangular GLT sequences. {\it  Electron. Trans. Numer. Anal.} \textbf{55} (2022), 585-617.

\bibitem{Barb 1D}   G. Barbarino, C. Garoni, S. Serra-Capizzano.  Block generalized locally Toeplitz sequences: theory and applications in the unidimensional case. {\it  Electron. Trans. Numer. Anal.} \textbf{53} (2020), 28-112.

\bibitem{Barb dD}
G. Barbarino, C. Garoni, S. Serra-Capizzano.  Block generalized locally Toeplitz sequences: theory and applications in the multidimensional case.  {\it  Electron. Trans. Numer. Anal.} \textbf{53} (2020), 113-216.

\bibitem{Barb NonH} G. Barbarino, S. Serra-Capizzano.  Non-Hermitian perturbations of Hermitian matrix-sequences and applications to the spectral analysis of the numerical approximation of partial differential equations. {\it  Numer. Linear Algebra Appl.} \textbf{27-3} (2020), e2286, 31 pp.

\bibitem{BeSe}
B. Beckermann, S. Serra-Capizzano.
On the asymptotic spectrum of Finite Element matrix sequences.
\SIAMJNA \textbf{45} (2007), 746-769.

\bibitem{Bene parallel}   P. Benedusi, P. Ferrari, C. Garoni, R. Krause,  S. Serra-Capizzano. Fast parallel solver for the space-time IgA-DG discretization of the diffusion equation. {\it  J. Sci. Comput.} \textbf{89-1} (2021), Paper No. 20, 21 pp.

\bibitem{Bene Rognes}
P. Benedusi, P. Ferrari, E. Rognes,  S. Serra-Capizzano. Modeling excitable cells with the EMI equations: spectral analysis and fast solution strategy. {\it  J. Sci. Comput.} \textbf{98-3} (2024), Paper No. 58, 23 pp.


\bibitem{Bene teo}   P. Benedusi, C. Garoni, R. Krause, X. Li,  S. Serra-Capizzano. Space-time FE-DG discretization of the anisotropic diffusion equation in any dimension: the spectral symbol. {\it  SIAM J. Matrix Anal. Appl.} \textbf{39-3} (2018), 1383-1420.


\bibitem{hidden} M. Benzi, D.A. Bini, D. Kressner, H. Munthe-Kaas, C. Van Loan.
\emph{Exploiting hidden structure in matrix computations: algorithms and applications.} Lecture notes from the CIME Summer School held in Cetraro, June 22–26, 2015. Edited by Michele Benzi and Valeria Simoncini. Lecture Notes in Mathematics, \textbf{2173}. Fondazione CIME/CIME Foundation Subseries. Springer, Cham; Centro Internazionale Matematico Estivo (C.I.M.E.), Florence, 2016.

\bibitem{Bhatia1997}
Bhatia, R.
\emph{Matrix Analysis.} Graduate Texts in Mathematics, Vol. 169, Springer-Verlag, New York, 1997.

\bibitem{bianchi}
D. Bianchi. Analysis of the spectral symbol associated to discretization schemes of linear self-adjoint differential operators. {\it  Calcolo} \textbf{58-3} (2021), Paper No. 38, 47 pp.

\bibitem{bianchi-S}
D. Bianchi, S. Serra-Capizzano. Spectral analysis of finite-dimensional approximations of 1d waves in non-uniform grids. {\it  Calcolo} \textbf{55-4} (2018), Paper No. 47, 28 pp.

\bibitem{Toe-use1} D.A. Bini, G. Latouche, B. Meini. \emph{Numerical methods for structured Markov chains.}
 Numerical Mathematics and Scientific Computation. Oxford Science Publications. Oxford University Press, New York, 2005.


\bibitem{BogoFDE1}
M. Bogoya, S. Grudsky, M. Mazza, S. Serra-Capizzano. On the extreme eigenvalues and asymptotic conditioning of a class of Toeplitz matrix-sequences arising from fractional problems. {\it  Linear Multilinear Algebra} {\bf 71-15}  (2023),  2462–2473.

\bibitem{BGS} M. Bogoya, S. Grudsky, S. Serra-Capizzano.
Fast non-Hermitian Toeplitz eigenvalue computations, joining matrixless algorithms and FDE approximation matrices.
 {\it SIAM J. Matrix Anal. Appl.} \textbf{45-1} (2024), 284-305.

\bibitem{BGS-2} M. Bogoya, S. Grudsky, S. Serra-Capizzano.
Eigenvalue superposition for Toeplitz matrix-sequences with matrix order dependent symbols.
 \LAA \textbf{697} (2024), 487-527.


\bibitem{BogoFDE2}
M. Bogoya, S. Grudsky, S. Serra-Capizzano, C. Tablino-Possio. Fine spectral estimates with applications to the optimally fast solution of large FDE linear systems. {\it BIT} {\bf 62-4}   (2022), 1417–1431.

\bibitem{mom-Toe}
M. Bolten, S.-E. Ekström, I. Furci, S. Serra-Capizzano.  Toeplitz momentary symbols: definition, results, and limitations in the spectral analysis of structured matrices. {\it Linear Algebra Appl.} \textbf{651}  (2022), 51–82.

\bibitem{mom-GLT}
M. Bolten, S.-E. Ekström, I. Furci, S. Serra-Capizzano.  A note on the spectral analysis of matrix sequences via GLT momentary symbols: from all-at-once solution of parabolic problems to distributed fractional order matrices. {\it Electron. Trans. Numer. Anal.} \textbf{58} (2023), 136–163. 

\bibitem{BS}
A. B\"ottcher, B. Silbermann.
\emph{Introduction to Large Truncated Toeplitz Matrices.}
Springer-Verlag, New York, 1999.


\bibitem{Ivrii} M. Bronstein, V. Ivrii. Sharp spectral asymptotics for operators with irregular coefficients. Pushing the limits, {\it Comm. Partial Differential Equations},  {\bf 28-1/2} (2003), 83-102.



\bibitem{T-Chan SIREV}
T.F. Chan, H.C. Elman. Fourier analysis of iterative methods for elliptic problems. {\it  SIAM Rev.} \textbf{31-1} (1989),  20-49.

\bibitem{Serra-koro-NM}
F. Di Benedetto, S. Serra-Capizzano. A unifying approach to abstract matrix algebra preconditioning. {\it
Numer. Math.} \textbf{82-1} (1999), 57-90.

\bibitem{Serra-koro-lama}
F. Di Benedetto, S. Serra-Capizzano.
Optimal multilevel matrix algebra operators. {\it  Linear Multilin. Algebra} \textbf{48-1} (2000), 35-66.

\bibitem{DGMSS}
M. Donatelli, C. Garoni, C. Manni, S. Serra-Capizzano, H. Speleers.
Robust and optimal multi-iterative techniques for IgA Galerkin linear systems.
\CMAME \textbf{284} (2015), 230-264.

\bibitem{DGMSScol}
M. Donatelli, C. Garoni, C. Manni, S. Serra-Capizzano, H. Speleers.
Robust and optimal multi-iterative techniques for IgA collocation linear systems.
\CMAME \textbf{284} (2015), 1120-1146.

\bibitem{our-MC}
M. Donatelli, C. Garoni, C. Manni, S. Serra-Capizzano, H. Speleers.
Spectral analysis and spectral symbol of matrices in isogeometric collocation methods.
\MC \textbf{85-300}  (2016),  1639-1680.

\bibitem{garoni-IgA-multigrid}
M. Donatelli, C. Garoni, C. Manni, S. Serra-Capizzano, H. Speleers.
 Symbol-based construction and analysis of optimal/robust multigrid methods for B-splines Isogeometric Analysis.
\SIAMJNA  \textbf{55-1}  (2017), 31-62.



\bibitem{fract-derivatives}
M. Donatelli,  M. Mazza, S. Serra-Capizzano.
 Spectral analysis and preconditioning for variable coefficient fractional derivative operators.
\JCP  \textbf{307}  (2016),  262-279.

\bibitem{fract-derivatives2} M. Donatelli, M. Mazza, S. Serra-Capizzano. Spectral analysis and multigrid methods for finite volume approximations of space-fractional diffusion equations. {\it SIAM J. Sci. Comput.} \textbf{40-6} (2018), A4007-A4039.


\bibitem{DNS} M. Donatelli, M. Neytcheva, S. Serra-Capizzano. {Canonical eigenvalue distribution of multilevel block Toeplitz sequences with non-Hermitian symbols.} Spectral theory, mathematical system theory, evolution equations, differential and difference equations, Birkhäuser/Springer Basel AG, Basel.
 {\it Oper. Theory Adv. Appl.} \textbf{221} (2012), 269–291.


\bibitem{DoroNS}
A. Dorostkar, M. Neytcheva, S. Serra-Capizzano.
Spectral analysis of coupled PDEs and of their Schur complements
via the notion of Generalized Locally Toeplitz sequences.
%
\CMAME \textbf{309} (2016),  74-105.


\bibitem{Ivo} I. Dravins, S. Serra-Capizzano, M. Neytcheva.
Spectral analysis of preconditioned matrices arising from stage-parallel implicit Runge-Kutta methods of arbitrarily high order.
 {\it SIAM J. Matrix Anal. Appl.} \textbf{45-2} (2024), 1007-1034.


\bibitem{NS DG Dumb} M. Dumbser, F. Fambri, I. Furci, M. Mazza, S. Serra-Capizzano, M. Tavelli. Staggered discontinuous Galerkin methods for the incompressible Navier-Stokes equations: spectral analysis and computational results. {\it Numer. Linear Algebra Appl.} \textbf{25-5} (2018), e2151, 31 pp.


\bibitem{eig expansion Laplacian} S.-E. Ekstr\"om, I. Furci, C. Garoni, C. Manni, S. Serra-Capizzano, H.  Speleers. Are the eigenvalues of the B-spline isogeometric analysis approximation of $-\Delta u=\lambda u$ known in almost closed form?  {\it Numer. Linear Algebra Appl.} \textbf{25-5} (2018), e2198, 34 pp.

\bibitem{sven-erik-asymp-exp1}
S.-E. Ekstr\"om, C. Garoni, S. Serra-Capizzano.
 Are the eigenvalues of banded symmetric Toeplitz matrices known  in close form?  {\it Exp. Math.} \textbf{4} (2018), 478-487.

\bibitem{dyad-new1}
S.-E. Ekström, S. Serra-Capizzano. Eigenvalue isogeometric approximations based on B-splines: tools and results. {\it Advanced methods for geometric modeling and numerical simulation}, 57–76, Springer INdAM Ser., \textbf{35}, Springer, Cham, 2019.

\bibitem{FT}
D. Fasino, P. Tilli. Spectral clustering properties of block multilevel Hankel matrices. 
\LAA {\bf 306-1/3} (2000), 155–163.

\bibitem{stat-conv1} H. Fast. Sur la convergence statistique. {\it Colloq. Math.}, {\bf 2} (1951), 241–244.


\bibitem{flipped1} P. Ferrari, I. Furci, S. Hon, M. Mursaleen, S. Serra-Capizzano. The eigenvalue distribution of special 2-by-2 block matrix-sequences with applications to the case of symmetrized Toeplitz structures. {\it SIAM J. Matrix Anal. Appl.} \textbf{40-3} (2019), 1066-1086.

\bibitem{flipped-BCs}
P. Ferrari, I. Furci, S. Serra-Capizzano.
Flipped structured matrix-sequences in image deblurring with reflective and anti-reflective boundary conditions.
{\it Numer. Alg.} (2025), https://doi.org/10.1007/s11075-024-01960-3.




\bibitem{Fle1} J. Fleckinger-Pellé. The vibrations of a drum with fractal boundary. In Mathematics and the 21st century (Cairo, 2000), 305-322, World Sci. Publishing, River Edge, NJ, 2001.

\bibitem{Fle2} J. Fleckinger, G. M\'etivier. Th\'eorie spectrale des op\'erateurs uniform\'ement elliptiques sur quelques ouverts irr\'eguliers.
{\it C.R. Acad. Sci. Paris S\'er. A} {\bf 276} (1973), 913-916.


\bibitem{blocking - prequel}
I. Furci, A. Adriani, S. Serra-Capizzano.
Block structured matrix-sequences and their spectral and singular value canonical distributions: a general theory.
 {\it arXiv} (2025), arXiv:2501.14874.	


\bibitem{topo-G} C. Garoni. Topological foundations of an asymptotic approximation theory for sequences of matrices with increasing size. \LAA  \textbf{513} (2017), 324–341.


\bibitem{GMPSS}
C. Garoni, C. Manni, F. Pelosi, S. Serra-Capizzano, H. Speleers.
On the spectrum of stiffness matrices arising from isogeometric analysis.
\NM \textbf{127} (2014), 751-799.

\bibitem{full-Galerkin}
C. Garoni, C. Manni, S. Serra-Capizzano, D. Sesana, H. Speleers.
Spectral analysis and spectral symbol of matrices in isogeometric Galerkin methods.
\MC \textbf{86-305}   (2017), 1343-1373.

\bibitem{lusin-garoni}
C. Garoni, C. Manni, S. Serra-Capizzano, D. Sesana, H. Speleers.
Lusin theorem, GLT sequences and matrix computations: an application to the spectral analysis of PDE discretization matrices.
 {\it J. Math. Anal. Appl.}  \textbf{446-1 } (2017),  365-382.


\bibitem{immersed NLAA-1}
C. Garoni, C. Manni, F. Pelosi, H. Speleers. {Spectral analysis of matrices resulting from isogeometric immersed methods and trimmed geometries.}
{\it  Comput. Methods Appl. Mech. Eng.} {\bf 400} (2022), paper 115551.

\bibitem{glt-garoni}
C. Garoni, S. Serra-Capizzano.
\emph{The theory of Generalized Locally Toeplitz sequences:  theory and applications - Vol I.}
SPRINGER - Springer Monographs in Mathematics,  Berlin, 2017.


\bibitem{glt-garoni-vol2}
C. Garoni, S. Serra-Capizzano.
\emph{The theory of Generalized Locally Toeplitz sequences: theory and applications - Vol II.}
SPRINGER - Springer Monographs in Mathematics,  Berlin, 2018.
%


\bibitem{glt-garoni-CIME}
C. Garoni, S. Serra-Capizzano.
\emph{Generalized Locally Toeplitz Sequences: A Spectral Analysis Tool for Discretized Differential Equations.}
  SPRINGER - Lecture Notes in Mathematics, CIME Foundation Subseries \textbf{2219} (2018), 161-236.


\bibitem{FEM}
C. Garoni, S. Serra-Capizzano, D. Sesana.
Spectral Analysis and Spectral Symbol of $d$-variate $Q_p$ Lagrangian FEM Stiffness Matrices.
\SIAMJMAA \textbf{36-3} (2015),  1100-1128

\bibitem{hughes}
C. Garoni, H. Speleers,  S.-E. Ekstr\"om, A. Reali, S. Serra-Capizzano, T.J.R. Hughes. Symbol-based analysis of finite element and isogeometric B-spline discretizations of eigenvalue problems: exposition and review. {\it
Arch. Comput. Methods Eng.} \textbf{26-5} (2019), 1639-1690.


\bibitem{Gil} P. Gilkey. {\em Asymptotic formulae in spectral geometry}. Studies in Advanced Mathematics, Chapman and Hall/CRC, Boca Raton, FL, 2004.


\bibitem{GoUMI} L. Golinskii, K. Kumar, M.N.N. Namboodiri, S. Serra-Capizzano. A note on a discrete version of Borg's theorem via Toeplitz-Laurent operators with matrix-valued symbols. {\it  Boll. Unione Mat. Ital.} \textbf{(9)6-1} (2013), 205–218. 

\bibitem{GoSe} L. Golinskii, S. Serra-Capizzano. The asymptotic properties of the spectrum of nonsymmetrically perturbed Jacobi matrix sequences. 
{\it J. Approx. Theory} \textbf{144-1}  (2007), 84–102.

\bibitem{answer:kac} C. Gordon, D. Webb, S. Wolpert. Isospectral plane domains and surfaces via Riemannian orbifolds.  {\it Invent. Math.}, {\bf 110-1} (1992), 1-22.

\bibitem{Gru} G. Grubb. {\em Functional calculus of pseudodifferential boundary problems}. Progress in Mathematics, Birkhäuser, Basel, 1996.

\bibitem{holm-serra} S. Holmgren, S. Serra-Capizzano, P. Sundqvist. Can one hear the composition of a drum? {\it Mediterr. J. Math.}, {\bf 3-2}  (2006), 227-249.


\bibitem{flipped PDE} S. Hon, S. Serra-Capizzano.  Block Toeplitz preconditioners for all-at-once systems from linear wave equations. {\it Electr. Trans. Numer.  Anal.} \textbf{58} (2023), 177-195.

\bibitem{flipped wave control} S. Hon, J. Dong, S. Serra-Capizzano. A preconditioned MINRES method for optimal control of wave equations and its related asymptotic spectral distribution theory. {\it SIAM J. Matrix Anal. Appl.} \textbf{44-4} (2023), 1477-1509.


\bibitem{hormander}
L. H\"ormander. Pseudo-differential operators and non-elliptic boundary problems. {\it Annals of Math.}, \textbf{2-83} (1966), 129-209.

\bibitem{GM-GLT}
A. Ilyas, M.F. Khan, V. Loi, S. Serra-Capizzano. A maximal distribution result and numerical tests for geometric means of HPD GLT matrix-sequences with degenerate commuting/non-commuting GLT symbols. \LAA (2025).

\bibitem{question:kac} M. Kac. Can one hear the shape of a drum? {\it Amer. Math. Monthly}, {\bf 73-4} (1966), 1-23.


\bibitem{koro-infinity-1} K. Kumar, M.N.N. Namboodiri, S. Serra-Capizzano. Perturbation of operators and approximation of spectrum. {\it
Proc. Indian Acad. Sci. Math. Sci.}  {\bf 124-2} (2014), 205–224.

\bibitem{koro-infinity-2} K. Kumar, M.N.N. Namboodiri, S. Serra-Capizzano. Preconditioners and Korovkin-type theorems for infinite-dimensional bounded linear operators via completely positive maps. {\it Studia Math.}   {\bf 218-2}  (2013), 95–118.


\bibitem{Kuij}
A.B.J. Kuijlaars, S. Serra-Capizzano. Asymptotic zero distribution of orthogonal polynomials with discontinuously varying recurrence coefficients. {\it J. Approx. Theory.} \textbf{113} (2001), 142-155.


\bibitem{stile-bianchi}
N. Lamsahel, A.E. Akri, A. Ratnani. Eigenvalues Distributions and Control Theory. {\it arXiv} (2024),  arXiv:2401.01975.

\bibitem{outl-free}
N. Lamsahel, C. Manni, A. Ratnani, S. Serra-Capizzano, H. Speleers. Outlier-free isogeometric discretizations for Laplace eigenvalue problems: closed-form eigenvalue and eigenvector expressions., \NM (2025), to appear.

\bibitem{toe-quaternion} X.-L. Lin, M.K. Ng, J. Pan. Hermitian Quaternion Toeplitz Matrices by Quaternion-valued Generating Functions. {\it arXiv} (2025), arXiv.2504.15073.

\bibitem{zuazua}
A. Marica, E. Zuazua. Propagation of 1D waves in regular discrete heterogeneous media: a Wigner measure approach. {\it Found. Comput. Math.} \textbf{15-6} (2015), 1571-1636.


\bibitem{curl-div IgA GLT} M. Mazza, C. Manni, A. Ratnani, S. Serra-Capizzano, H. Speleers. Isogeometric analysis for 2D and 3D curl-div problems: spectral symbols and fast iterative solvers. {\it Comput. Methods Appl. Mech. Engrg.} \textbf{344} (2019), 970-997.


\bibitem{flipping GLT1}
M. Mazza, J. Pestana. Spectral properties of flipped Toeplitz matrices and related preconditioning.
{\it BIT} \textbf{59-2} (2019), 463-482.


\bibitem{flipping GLT2}
M. Mazza, J. Pestana. The asymptotic spectrum of flipped multilevel Toeplitz matrices and of certain preconditionings.
{\it SIAM J. Matrix Anal. Appl.} \textbf{42-3} (2021), 1319-1336.


\bibitem{curl-curl IgA GLT} M. Mazza, A. Ratnani, S. Serra-Capizzano.  Spectral analysis and spectral symbol for the 2D curl-curl (stabilized) operator with applications to the related iterative solutions. {\it  Math. Comp.} \textbf{88-317} (2019), 1155-1188.

\bibitem{NS Schur GLT} M. Mazza, M. Semplice, S. Serra-Capizzano, E. Travaglia. A matrix-theoretic spectral analysis of incompressible Navier-Stokes staggered DG approximations and a related spectrally based preconditioning approach. {\it Numer. Math.} \textbf{149-4}  (2021),  933-971.

\bibitem{fract-derivatives3} M. Mazza, S. Serra-Capizzano, M. Usman. Symbol-based preconditioning for Riesz distributed-order space-fractional diffusion equations. {\it Electron. Trans. Numer. Anal.} \textbf{54} (2021), 499-513.



\bibitem{immersed NLAA-2}
M. Mazza, S. Serra-Capizzano, R.L. Sormani. { Spectral analysis and preconditioners for two-dimensional Riesz distributed-order space-fractional diffusion equations.} {\it  Numer. Linear Algebra Appl.} {\bf 31-3} (2023), e2536.


\bibitem{stat-conv4} M.A. Mursaleen, S. Serra-Capizzano. Statistical convergence via $q$-calculus and a Korovkin's type approximation theorem. {\it Axioms}, {\bf 11-2} (2022), Paper No. 70.

\bibitem{stat-conv3} M. Mursaleen, O.H.H. Edely. Statistical convergence of double sequences. {\it J. Math. Anal. Appl.}, {\bf 288} (2003), 223–231.

\bibitem{stat-conv5} M. Mursaleen, S.A. Mohiuddine. On lacunary statistical convergence with respect to the intuitionistic fuzzy normed space. {\it J. Comput. Appl. Math.}, {\bf 233-2} (2009), 142–149.

\bibitem{Toe-use2}
M.K. Ng. \emph{Iterative methods for Toeplitz systems.} Numerical Mathematics and Scientific Computation. Oxford University Press, New York, 2004.


\bibitem{Ngondiep} E. Ngondiep, S. Serra-Capizzano, D. Sesana. Spectral features and asymptotic properties for g-circulants and g-Toeplitz sequences. {\it  SIAM J. Matrix Anal. Appl.} \textbf{31-4} (2009/10), 1663-1687.

\bibitem{pre-LT}
S.V. Parter. On the eigenvalues of certain generalizations of Toeplitz matrices. {\it  Arch. Rat. Math.
Mech.} \textbf{3} (1962), 244–257.

\bibitem{flipped2} J. Pestana, A. Wathen. A preconditioned MINRES method for nonsymmetric Toeplitz matrices. {\it  SIAM J. Matrix Anal. Appl.} \textbf{36-1} (2015), 273-288.

\bibitem{SSS}
E. Salinelli, S. Serra-Capizzano, D. Sesana. Eigenvalue-eigenvector structure of Schoenmakers-Coffey matrices via Toeplitz technology and applications. 
\LAA {\bf 491} (2016), 138–160.



\bibitem{Sarathkumar}
N.S. Sarathkumar, S. Serra-Capizzano. GLT sequences and automatic computation of the symbol. (English summary)
\LAA \textbf{697}  (2024), 468–486.


\bibitem{Alec Beta} A.J.A. Schiavoni-Piazza, D. Meadon, S. Serra-Capizzano.
The $\beta$ maps: Strong clustering and distribution results on the complex unit circle.
\LAA \textbf{697} (2024), 365-383.

\bibitem{Alec cluster} A.J.A. Schiavoni-Piazza, S. Serra-Capizzano. Distribution results for a special class of matrix sequences: joining approximation theory and asymptotic linear algebra. {\it  Electron. Trans. Numer. Anal.} \textbf{59} (2023), 1-8.


\bibitem{See} R. Seeley. Complex powers of an elliptic operator. In  Singular Integrals (Proc. Sympos. Pure Math., Chicago, Ill., 1966), 288-307, AMS, Providence, RI.




\bibitem{Serra-multi-it}
S. Serra-Capizzano.
Multi-iterative methods. {\it  Comput. Math. Appl.} \textbf{26-4} (1993), 65-87.

\bibitem{Serra96}
S. Serra-Capizzano. Preconditioning strategies for Hermitian Toeplitz systems with nondefinite generating functions. {\it  SIAM J. Matrix Anal. Appl.} \textbf{17-4} (1996), 1007-1019.

\bibitem{Serra-taud1}
S. Serra-Capizzano.
The rate of convergence of Toeplitz based PCG methods for second order nonlinear boundary value problems. \NM \textbf{81-3} (1999), 461-495.


\bibitem{Serra-koro1}
S. Serra-Capizzano.
A Korovkin-type theory for finite Toeplitz operators via matrix algebras. \NM \textbf{82-1} (1999), 117-142.

\bibitem{Serra-koro2}
S. Serra-Capizzano.
A Korovkin-based approximation of multilevel Toeplitz matrices (with rectangular unstructured blocks) via multilevel trigonometric matrix spaces. {\it SIAM J. Numer. Anal.} \textbf{36-6} (1999), 1831–1857.

\bibitem{Serra99}
S. Serra-Capizzano. Spectral and computational analysis of block Toeplitz matrices having nonnegative definite matrix-valued generating functions. {\it  BIT} \textbf{39-1} (1999), 152-175.

\bibitem{Serra-acs}
S. Serra-Capizzano. Distribution results on the algebra generated by Toeplitz sequences: a finite-dimensional approach. \LAA \textbf{328-1/3}  (2001),  121-130.

\bibitem{Serra-taud2}
S. Serra-Capizzano.  Spectral behavior of matrix sequences and discretized boundary value problems. \LAA \textbf{337}  (2001), 37-78.


\bibitem{Serra03}
S. Serra-Capizzano.
Generalized Locally Toeplitz sequences: spectral analysis and applications to discretized partial differential equations.
\LAA \textbf{366} (2003), 371-402.

\bibitem{Serra06}
S. Serra-Capizzano.
The GLT class as a generalized Fourier Analysis and applications.
\LAA \textbf{419} (2006), 180-233.

\bibitem{dyad-1}
S. Serra-Capizzano, C. Tablino-Possio. Cristina Spectral and structural analysis of high precision finite difference matrices for elliptic operators. {\it Linear Algebra Appl.} \textbf{293-1/3} (1999), 85–131.

\bibitem{dyad-2}
S. Serra-Capizzano, C. Tablino-Possio.  Positive representation formulas for finite difference discretizations of (elliptic) second order PDEs. {\it Structured matrices in mathematics, computer science, and engineering, II} (Boulder, CO, 1999), 295–317, {\it  Contemp. Math.}, \textbf{281}, Amer. Math. Soc., Providence, RI, 2001.

\bibitem{dyad-3}
S. Serra-Capizzano, C. Tablino-Possio.  Preconditioning strategies for 2D finite difference matrix sequences. {\it Electron. Trans. Numer. Anal.} \textbf{16}  (2003), 1–29.


\bibitem{Serra-Tilli1}
S. Serra-Capizzano, P. Tilli. Extreme singular values and eigenvalues of non-Hermitian block Toeplitz matrices. {\it J. Comput. Appl. Math.} \textbf{108-1/2} (1999), 113-130.

\bibitem{Serra-Tilli2}
S. Serra-Capizzano, P. Tilli. On unitarily invariant norms of matrix-valued linear positive operators. {\it J. Inequal. Appl.} \textbf{7-3} (2002), 309-330.

\bibitem{nega-equi-u-trasf-algebras}
S. Serra-Capizzano, E.E. Tyrtyshnikov. How to prove that a preconditioner cannot be superlinear. \MC \textbf{72-243} (2003), 1305-1316.


\bibitem{TaoVu-RandMatr}
T. Tao, V. Vu. Random matrices: the circular law. {\it Commun. Contemp. Math.} \textbf{10-2} (2008), 261-307.

\bibitem{Tilli-ToeL1}
P. Tilli. A note on the spectral distribution of Toeplitz matrices. {\it Linear Multilin. Algebra} \textbf{45-2/3} (1998), 147-159.

\bibitem{Tilli}
P. Tilli.
Locally Toeplitz sequences: spectral properties and applications.
\LAA \textbf{278} (1998), 91-120.

\bibitem{Ty96}
E.E. Tyrtyshnikov.
A unifying approach to some old and new theorems on distribution and clustering.
\LAA {\bf 232} (1996), 1-43.

\bibitem{draft-T}
E.E. Tyrtyshnikov.
Extra dimension approach to Spectral Distributions. {\it Private discussion}, 1997.

\bibitem{ZamTy-ToeL1}
N.L. Zamarashkin, E.E.  Tyrtyshnikov. Distribution of the eigenvalues and singular numbers of Toeplitz matrices under weakened requirements on the generating function. {\it Mat. Sb.} {\bf 188-8} (1997), 83–92; English translation in {\it Sb. Mat.} {\bf 188-8} (1997), 1191-1201.

\end{thebibliography}
		\end{document}